\documentclass[11pt]{article}

\usepackage[utf8]{inputenc}
\usepackage{amsmath,amssymb,amsthm}
\usepackage{verbatim,xcolor,algorithm,algorithmic,xfrac}
\usepackage{graphicx}

\providecommand{\keywords}[1]
{
  \small	
  \textbf{\textit{Keywords---}} #1
}

\providecommand{\AMS}[1]
{
  \small	
  \textbf{\textit{AMS Subject Codes---}} #1
}

\newtheorem{theorem}{Theorem}
\newtheorem{lemma}{Lemma}
\newtheorem{remark}{Remark}
\newtheorem{proposition}{Proposition}
\newcommand{\R}{\mathbb{R}}
\newcommand{\tbeta}{\widetilde{\beta}}
\newcommand{\ty}{\widetilde{y}}
\newcommand{\tv}{\widetilde{v}}
\newcommand{\tV}{\widetilde{V}}
\newcommand{\tH}{\widetilde{H}}

\newcommand{\xkI}{x_k}
\newcommand{\xkW}{\tx_k}
\newcommand{\tx}{\widetilde{x}}
\newcommand{\rkI}{r_k}
\newcommand{\rkW}{\tr_k}
\newcommand{\tr}{\widetilde{r}}
\newcommand{\cE}{{\cal E}}
\newcommand{\be}{\begin{equation}}
\newcommand{\ee}{\end{equation}}
\newcommand{\bmx}{\begin{bmatrix}}
\newcommand{\emx}{\end{bmatrix}}

\title{Exploiting variable precision in GMRES%
 \thanks{
   The work of Ehouarn Simon was partially supported by
   the French National program LEFE
	 (Les Enveloppes Fluides et l’Environnement).
	 The work of Serge Gratton and Phillipe Toint
	 was partially supported by the 3IA Artificial and Natural Intelligence
	 Toulouse Institute, French ``Investing for the Future - PIA3'' program
	 under the Grant agreement ANR-19-PI3A-0004.}}

\author{%
   Serge Gratton\footnotemark[1] \and
   Ehouarn Simon\footnotemark[2] \and
   David Titley-Peloquin\footnotemark[3] \and
   Philippe Toint\footnotemark[4]
}

\begin{document}
\maketitle

\renewcommand{\thefootnote}{\fnsymbol{footnote}}

\footnotetext[1]{INPT-IRIT-ENSEEIHT, Toulouse, France (serge.gratton@toulouse-inp.fr).}
\footnotetext[2]{INPT-IRIT-ENSEEIHT, Toulouse, France (ehouarn.simon@toulouse-inp.fr).}
\footnotetext[3]{McGill University, Montreal, Canada (david.titley-peloquin@mcgill.ca).}
\footnotetext[4]{The University of Namur, Namur, Belgium (philippe.toint@unamur.be).}

\renewcommand{\thefootnote}{\arabic{footnote}}

\begin{abstract}
We describe how variable precision floating-point arithmetic
can be used to compute inner products in the iterative solver GMRES.
We show how the precision of the inner products
carried out in the algorithm can be reduced as the iterations proceed,
without affecting the convergence rate
or final accuracy achieved by the iterates.
Our analysis explicitly takes into account
the resulting loss of orthogonality in the Arnoldi vectors.
We also show how inexact matrix-vector products 
can be incorporated into this setting.
\end{abstract}


\AMS{
15A06, 65F10, 65F25, 97N20}

\pagestyle{myheadings}
\thispagestyle{plain}
\markboth{S. GRATTON, E. SIMON, D. TITLEY-PELOQUIN, AND P. TOINT}{INEXACT GMRES}

\section{Introduction}

As highlighted in a recent SIAM News article~\cite{High17},
there is growing interest in the use of variable precision
floating-point arithmetic in numerical algorithms.
(Other recent references
include~\cite{CarsHigh18,haid2018b,haid2018,haid2017,HighPran19b,HighPran19}
to cite only a few.)
In this paper, we describe how variable precision arithmetic
can be exploited in the iterative solver GMRES~\cite{SaadSchu86}.
We show that the precision of some floating-point operations
carried out in the algorithm can be reduced as the iterations proceed,
without affecting the convergence rate
or final accuracy achieved by the iterates.

There is already a literature on the use of inexact matrix-vector products
in GMRES and other Krylov subspace methods; see,
e.g.,~\cite{SimoSzyl03,VDEhSlei04,BourFray06,GiraGratLang07,SimoSzyl07,GratSimoToin18}
and the references therein.
This work is not a simple extension of such results.
To illustrate, suppose that all arithmetic operations are performed exactly,
except the matrix-vector products.
Then one obtains an inexact Arnoldi relation
\be\label{eq:1}
A V_k + E_k = V_{k+1} H_k, \qquad V_k^T V_k = I.
\ee
On the other hand, if only inner products are performed inexactly,
the Arnoldi relation continues to hold 
but the orthogonality of the Arnoldi vectors is lost:
\be\label{eq:2}
A V_k = V_{k+1} H_k, \qquad V_k^T V_k = I - F_k.
\ee
Thus, to understand the convergence behaviour
and maximum attainable accuracy of GMRES
implemented with inexact inner products,
it is absolutely necessary to understand
the resulting loss of orthogonality in the Arnoldi vectors.
We adapt techniques used in the rounding-error analysis
of the Modified Gram-Schmidt (MGS) algorithm
(see~\cite{Bjor67,BjorPaig92} 
or~\cite{LeonBjorGand13} for a more recent survey)
and of the MGS-GMRES algorithm
(see~\cite{Drko95,Gree97,PaigRozlStra06}).

We focus on inexact inner products and matrix-vector products
(as opposed to the other saxpy operations involved in the algorithm)
because these are the two most time-consuming operations
in parallel computations.
The rest of the paper is organized as follows.  
We start with a brief discussion of GMRES in non-standard inner products
in Section~\ref{sec:gmres}.
Next, in Section~\ref{sec:IP}, we analyze GMRES with inexact inner products.
We then show how inexact matrix-vector products
can be incorporated into this setting in Section~\ref{sec:matvec}.
Some numerical examples are presented in Sections~\ref{sec:num}
and~\ref{sec:mult}.

\section{GMRES in weighted inner products}\label{sec:gmres}

Shown below is the Arnoldi algorithm,
with $\langle y , z \rangle = y^Tz$
denoting the standard Euclidean inner product.

\begin{algorithm}
\caption{Arnoldi algorithm}
\label{alg:arnoldi}
\begin{algorithmic}[1]
\REQUIRE{$A\in\R^{n\times n}$, $b\in\R^n$}
\STATE $\beta = \sqrt{ \langle b , b \rangle }$
\STATE $v_1 = b/\beta$
\FOR{$j=1,2,\dots$}
\STATE $w_{j} = Av_j$
\FOR{$i=1,\dots,j$}
\STATE $h_{ij} = \langle v_i, w_j \rangle$
\STATE $w_{j} = w_{j} - h_{ij}v_i$
\ENDFOR
\STATE $h_{j+1,j} = \sqrt{ \langle w_{j} , w_{j} \rangle }$
\STATE $v_{j+1} = w_{j} / h_{j+1,j}$
\ENDFOR
\end{algorithmic}
\end{algorithm}

After $k$ steps of the algorithm are performed
in exact arithmetic, the output is 
$V_{k+1}=[v_1,\dots,v_{k+1}]\in\R^{n\times(k+1)}$
and upper-Hessenberg $H_k\in\R^{(k+1)\times k}$
such that 
$$
v_1 = \frac{b}{\beta},
  \qquad AV_k = V_{k+1}H_k,
  \qquad V_k^T V_k = I_k.
$$
The columns of $V_k$ form an orthonormal basis
for the Krylov subspace
$$
{\cal K}_k(A,b) = \mathrm{span}\big\{\, b,Ab,A^2b,\dots,A^{k-1}b\,\big\}.
$$
In GMRES, we restrict $x_k$ to this subspace:
$x_k=V_ky_k$, where $y_k\in\R^k$ is the solution of
$$
\min_y \|b-AV_ky\|_{2}
   = \min_y \|V_{k+1} (\beta e_1 - H_k y)\|_2
   = \min_y \|\beta e_1 - H_k y\|_2.  
$$
It follows that
\be\label{eq:xr}
\begin{split}
& x_k = V_ky_k = V_k (H_k^TH_k)^{-1}H_k^T (\beta e_1) 
      = V_k H_k^\dag (\beta e_1), \\
& r_k = b-Ax_k = V_{k+1} (\beta e_1 - H_k y_k) 
     = V_{k+1} (I-H_kH_k^\dag) \beta e_1.
\end{split}
\ee

Any given symmetric positive definite matrix $W$ defines
a weighted inner product $\langle y,z \rangle_W=y^TWz$
and associated norm $\|z\|_W = \sqrt{ \langle z,z \rangle_W }$. 
Suppose we use this inner product
instead of the standard Euclidean inner product in the Arnoldi algorithm.
We use tildes to denote the resulting quantities in the algorithm.
After $k$ steps, the result is $\tV_{k+1}=[\tv_1,\dots,\tv_{k+1}]$
and upper-Hessenberg $\tH_k\in\R^{(k+1)\times k}$
such that 
$$
\tv_1 = \frac{b}{\|b\|_{W}}
      = \frac{b}{\tbeta},
  \qquad A\tV_k = \tV_{k+1}\tH_k,
  \qquad \tV_k^T W \tV_k = I_k.
$$
The columns of $\tV_k$ form a $W$-orthonormal basis
for ${\cal K}_k(A,b)$.
Let $\tx_k=\tV_k\ty_k$,
where $\ty_k\in\R^k$ is the solution of
$$
\min_y \|b-A\tV_ky\|_{W}
   = \min_y \|\tV_{k+1} (\tbeta e_1 - \tH_k y)\|_{W}
   = \min_y \|\tbeta e_1 - \tH_k y\|_2,  
$$
so that 
$$
\tx_k = \tV_k\tH_k^\dag (\tbeta e_1), \qquad
\tr_k = b-A\tx_k = \tV_{k+1} (I-\tH_k\tH_k^\dag) \tbeta e_1.
$$
We denote the above algorithm $W$-GMRES.

Let $\xkI$ and $\xkW$ denote the iterates computed
by standard GMRES and $W$-GMRES, respectively,
with corresponding residual vectors $\rkI$ and $\rkW$.
It is well known that 
\be\label{eq:res}
1 \leq \frac{\|\rkW\|_2}{\|\rkI\|_2}
  \leq \sqrt{\kappa_2(W)}.
\ee
See e.g.~\cite{PestWath13} for a proof.
Thus, if $\kappa_2(W)$ is small, 
the Euclidean norm of the residual vector in $W$-GMRES converges
at essentially the same rate as in standard GMRES.
A similar result~\cite[Theorem 4]{Fre92b} holds
for the residual computed in
the quasi-minimal residual method~\cite{FreNac91,Fre92a}. 


\section{GMRES with inexact inner products}\label{sec:IP}

\subsection{Preliminary results}

Suppose  the inner products in the Arnoldi algorithm
are computed inexactly, i.e.,
line~6 in Algorithm~\ref{alg:arnoldi} is replaced by
\be\label{eq:eta}
h_{ij} = v_i^T w_{j} + \eta_{ij},
\ee
with $|\eta_{ij}|$ bounded by some tolerance.
Our main contribution is to show precisely how large each $\eta_{ij}$
can be without affecting the convergence of GMRES.
Throughout we assume that all arithmetic operations in GMRES
are performed exactly, except for the above inner products.

It is straightforward to show that
despite the inexact inner products in~\eqref{eq:eta},
the relation $AV_k = V_{k+1} H_k$ continues to hold.
On the other hand, the orthogonality
of the Arnoldi vectors is lost.
We have
\be\label{eq:IF}
[b,AV_k] = V_{k+1}[\beta e_1, H_k], \qquad
V_{k+1}^T V_{k+1} = I_{k+1} + F_k.
\ee
The relation between each $\eta_{ij}$
and the overall loss of orthogonality $F_k$
is very difficult to understand.
To simplify the analysis we suppose
that each $v_j$ is normalized exactly.
(This is not an uncommon assumption;
see, e.g.,~\cite{Bjor67} and~\cite{Paig09}.)
Under this simplification, we have 
\be\label{eq:F}
F_k = \bar{U}_k + \bar{U}_k^T, \quad
\bar{U}_k = \bmx 0_{k\times1} & U_k \\ 
                 0_{1\times1} & 0_{1\times k} \emx, \quad
U_k = \bmx v_1^Tv_2 & \dots  & v_1^Tv_{k+1} \\
                    & \ddots & \vdots       \\
                    &        & v_k^Tv_{k+1} \emx,
\ee
i.e., $U_k\in\R^{k\times k}$ contains
the strictly upper-triangular part of $F_k$.
Define
\be\label{eq:NR}
N_k =
  \bmx \eta_{11} & \dots  & \eta_{1k} \\
                 & \ddots & \vdots    \\
                 &        & \eta_{kk}
  \emx, \quad
R_k = 
  \bmx\, h_{21} & \dots  & h_{2k} \\
                & \ddots & \vdots    \\
                &        & h_{k+1,k}
  \emx.
\ee
Note that $R_k$ must be invertible if $h_{j+1,j}\neq 0$ for $j=1,\dots,k$,
in other words, if GMRES has not terminated by step $k$.
(We assume that GMRES does not breakdown by step $k$.)
Following Bj\"orck's seminal rounding error analysis
of MGS~\cite{Bjor67}, it can be shown that
\be\label{eq:NUR}
N_k = - [ 0, U_k ] H_k = - U_k R_k.
\ee
For completeness, a proof of~\eqref{eq:NUR} is provided in the appendix.

Additionally, in order to understand how $\|F_k\|_2$ increases
as the residual norm decreases, 
we will need the following rather technical lemma.
The relationship~\eqref{eq:ejyk} is well know
(see for example~\cite[Lemma~5.1]{SimoSzyl03})
while~\eqref{eq:STLS} is essentially
a special case of~\cite[Theorem~4.1]{PaigStra02a}.
We defer the proof to the appendix.

\medskip
\begin{lemma}\label{lem:STLS}
Let $y_j$ and $t_j$ be the least squares solution
and residual vector of
$$
\min_{y} \|\beta e_1 - H_jy\|_2,
$$
for $j=1,\dots,k$.
Then 
\be\label{eq:ejyk}
|e_j^T y_k| \leq \frac{\|t_{j-1}\|_2}{\sigma_{\min}(H_k)}.
\ee
In addition, given $\epsilon>0$, let $D_k$ be any nonsingular matrix such that 
\be\label{eq:condD}
\|D_k\|_2 \leq  \frac{ \sigma_{\min}(H_k) \epsilon \|b\|_2 }{\sqrt{2}\|t_k\|_2}.
\ee
Then
\be\label{eq:STLS}
\frac{\|t_k\|_2}{ \big( \epsilon^2\|b\|_2^2 + 2\|D_ky_k\|_2^2 \big)^{1/2} }
  \leq \sigma_{\min} \left( \bmx \epsilon^{-1} e_1, H_kD_k^{-1} \emx \right)
  \leq \frac{\|t_k\|_2}{\epsilon\|b\|_2}.
\ee
\end{lemma}
\medskip

Finally, although the columns of $V_{k+1}$ in~\eqref{eq:IF}
are not orthonormal in the standard Euclidean inner product,
we will use the fact that there exists an inner product
in which they are orthonormal.
The proof of the following lemma is given in the appendix.

\medskip
\begin{lemma}\label{lem:qre}
Consider a given matrix $Q\in\R^{n\times k}$ of rank $k$ such that
\be\label{eq:q1}
\qquad Q^TQ = I_k - F.
\ee
If $\|F\|_2\leq\delta$ for some $\delta\in(0,1)$,
then there exists a matrix $M$
such that $I_n+M$ is symmetric positive definite and 
\be\label{eq:q2}
Q^T(I_n+M)Q = I_k. 
\ee
In other words, the columns of $Q$ are exactly
orthonormal in an inner product defined by $I_n+M$.
Furthermore,
\be\label{eq:kapdelta}
\kappa_2(I_n+M) \leq \frac{1+\delta}{1-\delta}.
\ee
\end{lemma}
\medskip

Note that $\kappa_2(I_n+M)$ remains small even 
for values of $\delta$ close to $1$.
For example, suppose $\|I_k-Q^TQ\|_2 = \delta = \sfrac{1}{2}$, 
indicating an extremely severe loss of orthogonality.
Then $\kappa_2(I_n+M)\leq 3$,
so $Q$ still has exactly orthonormal 
columns in an inner product defined
by a very well-conditioned matrix.

\medskip
\begin{remark}\label{rem:1}
Paige and his coauthors~\cite{BjorPaig92,Paig09,PaigWull14}
have developed an alternative measure of loss of orthogonality.
Given $Q\in\R^{n\times k}$ with normalized columns,
the measure is~$\|S\|_2$, where $S=(I+U)^{-1}U$
and $U$ is the strictly upper-triangular part of~$Q^TQ$.
Additionally, orthogonality can be recovered by augmentation: the matrix
$P=\left[\begin{smallmatrix} S \\ Q(I-S) \end{smallmatrix}\right]$
has orthonormal columns.
This measure was used in the groundbreaking rounding error analysis
of the MGS-GMRES algorithm~\cite{PaigRozlStra06}.
In the present paper, under the condition $\|F\|_2\leq\delta<1$,
we use the measure $\|F\|_2$ and recover orthogonality
in the $(I+M)$ inner product.
However, Paige's approach is likely to be the most appropriate
for analyzing the Lanczos and conjugate gradient algorithms,
in which orthogonality is quickly lost 
and $\|F\|_2>1$ long before convergence.
\end{remark}

\subsection{A strategy for bounding the $\eta_{ij}$}

We now show how to bound the error $\eta_{ij}$ in~\eqref{eq:eta} 
to ensure that the convergence of the GMRES 
is not affected by the inexact inner products.

The following theorem shows how the convergence of GMRES
with inexact inner products relates to that of exact GMRES.
The idea is similar to~\cite[Section 5]{PaigRozlStra06},
in which the quantity $\|E_k R_k^{-1}\|_F$ must be bounded,
where $R_k$ is the matrix in~\eqref{eq:NR}
and $E_k$ is a matrix containing rounding errors.

\medskip
\begin{theorem}\label{thm:picketa}
Let $x_k^{\mathrm{(e)}}$ denote the $k$-th iterate of standard GMRES,
performed exactly, with residual vector $r_k^{\mathrm{(e)}}$.
Now suppose that the Arnoldi algorithm is run with inexact inner products
as in~\eqref{eq:eta}, so that~\eqref{eq:IF}--\eqref{eq:NUR} hold,
and let $x_k$ and $r_k$ denote the resulting GMRES iterate and residual vector.
Let $y_k$ and $t_k$ be the least squares solution
and residual vector of
$$
\min_{y} \|\beta e_1 - H_ky\|_2.
$$
If for all steps $j=1,\dots,k$ of GMRES
all inner products are performed inexactly 
as in~\eqref{eq:eta} with tolerances bounded by
\be\label{eq:picketa}
|\eta_{ij}| \leq \eta_j 
  \equiv \frac{ \phi_j\epsilon \sigma_{\min}(H_k) }{\sqrt{2}}
         \frac{\|b\|_2}{\|t_{j-1}\|_2} 
\ee
for any $\epsilon\in(0,1)$ and any positive numbers $\phi_j$
such that $\sum_{j=1}^k \phi_j^2 \leq 1$,
then at step $k$ either
\be\label{eq:ip}
1 \leq \frac{\|r_k\|_2}{\|r_k^{\mathrm{(e)}}\|_2}
  \leq \sqrt{3},
\ee
or
\be\label{eq:converged}
\frac{\|t_k\|_2}{\|b\|_2} \leq 6k \epsilon,
\ee
implying that GMRES has converged
to a relative residual of $6k\epsilon$.
\end{theorem}
\medskip
\begin{proof}
If~\eqref{eq:picketa} holds, then in~\eqref{eq:NR}
$$
|N_k| \leq
  \bmx \eta_{1} & \eta_{2} & \dots  & \eta_{k} \\
                & \eta_{2} & \dots  & \eta_{k} \\
                &          & \ddots & \vdots   \\
                &          &        & \eta_{k}
  \emx
   = E_kD_k,
$$
where $E_k$ is an upper-triangular matrix 
containing only ones in its upper-triangular part,
so that $\|E_k\|_2 \leq k$, 
and $D_k=\mbox{diag}(\eta_1,\dots,\eta_k)$.
Then, 
\be\label{eq:bound2}
\begin{split}
\|N_kR_k^{-1}\|_2  
    &\leq  \|N_kD_k^{-1}\|_2 \|D_kR_k^{-1}\|_2 \\
    &\leq  \|E_k\|_2 \|D_kR_k^{-1}\|_2
    \leq  k \|(R_kD_k^{-1})^{-1}\|_2.
\end{split}
\ee
Let $h_k^T$ denote the first row of $H_k$,
so that $H_k = \left[ \begin{smallmatrix} h_k^T \\ R_k \end{smallmatrix}\right]$.
For any $\epsilon>0$ we have
\begin{align*}
\sigma_{\min}(R_kD_k^{-1})
 &= \, \min_{\|u\|_2=\|v\|_2=1} u^T R_kD_k^{-1} v \\
 &= \, \min_{\|u\|_2=\|v\|_2=1} \ [0,u^T] \bmx \epsilon^{-1} & h_k^T D_k^{-1} \\ 0 & R_kD_k^{-1} \emx
                               \bmx 0 \\ v \emx \\
 &\geq \, \min_{\|u\|_2=\|v\|_2=1} \
        u^T \bmx \epsilon^{-1} & h_k^T D_k^{-1} \\ 0 & R_kD_k^{-1} \emx v \\
 &= \, \sigma_{\min} \left( \bmx \epsilon^{-1} e_1, H_kD_k^{-1} \emx \right).
\end{align*}
Therefore,
$$
\| (R_kD_k^{-1})^{-1} \|_2
  = \frac{1}{\sigma_{\min}(R_kD_k^{-1})} \\
  \leq \frac{1}{\sigma_{\min}
       \left( \bmx \epsilon^{-1} e_1, H_kD_k^{-1} \emx \right)}.
$$
Notice that if the $\eta_j$ are chosen as in~\eqref{eq:picketa},
$D_k$ automatically satisfies~\eqref{eq:condD}.
Using the lower bound in~\eqref{eq:STLS},
then~\eqref{eq:ejyk} and~\eqref{eq:picketa}, we obtain
\begin{align*}
\| (R_kD_k^{-1})^{-1} \|_2
  &\leq \frac{ \big( \epsilon^2\|b\|_2^2 
         + 2 \|D_ky_k\|_2^2 \big)^{1/2} }{\|t_k\|_2} \\
  &= \frac{ \big( \epsilon^2\|b\|_2^2 
         + 2 \sum_{j=1}^k \eta_j^2 (e_j^Ty_k)^2 \big)^{1/2} }{\|t_k\|_2} \\
  &\leq \frac{ \big( \epsilon^2\|b\|_2^2 
       + \sum_{j=1}^k \phi_j^2\epsilon^2\|b\|_2^2 \big)^{1/2} }{\|t_k\|_2}
  = \frac{\sqrt{2}\epsilon\|b\|_2}{\|t_k\|_2}.
\end{align*}
Therefore, in~\eqref{eq:bound2},
$$
\|N_kR_k^{-1}\|_2 \leq \frac{\sqrt{2}k\epsilon\|b\|_2}{\|t_k\|_2}
                  \leq \frac{6k\epsilon\|b\|_2}{\|t_k\|_2} \frac{1}{4}.
$$
\medskip
If~\eqref{eq:converged} does not hold,
then $\|N_k^{\vphantom 2}R_k^{-1}\|_2 \leq \sfrac{1}{4}$.
From~\eqref{eq:F} and~\eqref{eq:NUR},
we have
\be\label{eq:bound}
\|F_k\|_2 \leq 2 \|U_k\|_2 = 2\|N_k R_k^{-1}\|_2,
\ee
with the matrix $F_k$ defined in~\eqref{eq:IF}.
Thus,
$\|F_k\|_2\leq \frac{1}{2} < 1$ and we can apply 
Lemma~\ref{lem:qre} with $Q = V_{k+1}$ and $\delta = \frac{1}{2}$.
There is a symmetric positive definite matrix~$W=I_n+M$ such that
$$
[b,AV_k] = V_{k+1} [\beta e_1, H_k], \quad
V_{k+1}^T W V_{k+1} = I_{k+1}, \quad
\kappa_2(W) \leq \frac{1+\delta}{1-\delta}=3.
$$
The Arnoldi algorithm implemented with inexact inner products
has computed an $W$-orthonormal basis for the Krylov subspace ${\cal K}_k(A,b)$.
The iterate $x_k$ is the same as the iterate that would have been obtained
by running $W$-GMRES exactly, and~\eqref{eq:res} implies~\eqref{eq:ip}.

Therefore, if the $|\eta_{ij}|$ are bounded by 
tolerances $\eta_j$ chosen as in~\eqref{eq:picketa},
either~\eqref{eq:ip} holds,
or~\eqref{eq:converged} holds.
\end{proof}
\medskip

Theorem~\ref{thm:picketa} can be interpreted as follows.
If at all steps $j=1,2,\dots$ of GMRES the inner products
are computed inaccurately with tolerances $\eta_j$ in~\eqref{eq:picketa},
then convergence at the same rate as exact GMRES
is achieved until a relative residual of essentially $k\epsilon$ is reached.
Notice that $\eta_j$ is inversely proportional to the residual norm.
This allows the inner products to be computed more and more
inaccurately as as the iterations proceed.

\subsection{Practical considerations}

If no more than $K_{\max}$ iterations are to be performed,
we can let $\phi_j=K_{\max}^{-\sfrac{1}{2}}$
(although more elaborate choices for $\phi_j$ could be considered;
see for example~\cite{GratSimoToin18}).
Then the factor $\sfrac{\phi_j}{\sqrt{2}}$ in~\eqref{eq:picketa}
can be absorbed along with the $k$ in~\eqref{eq:converged}.

One important difficulty with~\eqref{eq:picketa} is that $\sigma_{\min}(H_k)$
is required to pick $\eta_j$ at the start of step~$j$,
but $H_k$ is not available until the final step~$k$.
A similar problem occurs in GMRES with inexact matrix-vector
products; see~\cite{SimoSzyl03,VDEhSlei04}
and the comments in Section~\ref{sec:matvec}.
In our experience, is often possible
to replace $\sigma_{\min}(H_k)$ in~\eqref{eq:picketa} by $1$,
without significantly affecting the convergence of GMRES.  
This leads to following:
\be\label{eq:picketa2}
 \mbox{Aggressive threshhold}: 
     \qquad\quad \eta_j = \epsilon \frac{\|b\|_2}{\|t_{j-1}\|_2}, 
     \quad j=1,2,\dots.\qquad\quad\
\ee
In exact arithmetic, $\sigma_{\min}(H_k)$ is bounded below by $\sigma_{\min}(A)$. 
If the smallest singular value of $A$ is known, one can estimate
$\sigma_{\min}(H_k)\approx\sigma_{\min}(A)$ in~\eqref{eq:picketa},
leading to the following:
\be\label{eq:picketa3}
 \mbox{Conservative threshhold}:
    \qquad \eta_j = \epsilon\,\sigma_{\min}(A) \frac{\|b\|_2}{\|t_{j-1}\|_2}, 
    \quad j=1,2,\dots.
\ee
This prevents potential early stagnation of the residual norm,
but is often unnecessarily stringent.
(It goes without saying that if the conservative threshold 
is less than~$u\|A\|_2$, where $u$ is the machine precision,
then the criterion is vacuous:
according to this criterion no inexact inner products
can be carried out at iteration~$j$.)
Numerical examples are given in Sections~\ref{sec:num} and~\ref{sec:mult}.


\section{Incorporating inexact matrix-vector products}\label{sec:matvec}

As mentioned in the introduction, there is already a literature 
on the use of inexact matrix-vector products in GMRES.
These results are obtained by assuming
that the Arnoldi vectors are orthonormal
and analyzing the inexact Arnoldi relation
$$
A V_k + E_k = V_{k+1} H_k, \qquad V_k^T V_k = I.
$$
In practice, however, the computed Arnoldi vectors
are very far from being orthonormal, 
even when all computations are performed
in double precision arithmetic;
see for example~\cite{Drko95,Gree97,PaigRozlStra06}.

The purpose of this section is to show that
the framework used in~\cite{SimoSzyl03} and~\cite{VDEhSlei04}
to analyze inexact matrix-vector products in GMRES is still valid
when the orthogonality of the Arnoldi vectors is lost,
i.e., under the inexact Arnoldi relation
\be\label{eq:3}
A V_k + E_k = V_{k+1} H_k, \qquad V_k^T V_k = I - F_k.
\ee
We assume that the errors $\eta_{ij}$ in computing the inner products
is sufficiently small that $\|F_k\|_2\leq\delta < 1$,
as per Section~\ref{sec:IP}.
Then from Lemma~\ref{lem:qre} there exists
a symmetric positive definite matrix $W=I_n+M\in\R^{n\times n}$
such that $V_{k+1}^T W V_{k+1} = I_{k+1}$,
and with singular values bounded as in~\eqref{eq:kap}.

\subsection{Bounding the residual gap}

As in previous sections,
we use $x_k=V_ky_k$ to denote the computed GMRES iterate,
with $r_k=b-Ax_k$ for the actual residual vector
and $t_k = \beta_1 e_1 - H_ky_k$
for the residual vector updated in the GMRES iterations.
From 
$$
\|r_k\|_2 \leq \|r_k - V_{k+1}t_k\|_2 + \|V_{k+1}t_k\|_2,
$$
if 
\be\label{eq:gap1}
\max \left\{ \, \|r_k - V_{k+1}t_k\|_2, \, \|V_{k+1}t_k\|_2 \, \right\}
 \leq \frac{\epsilon}{2} \|b\|_2
\ee
then
\be\label{eq:gap2}
\|r_k\|_2 \leq \epsilon \|b\|_2.
\ee
From the fact that the columns of $W^{\sfrac{1}{2}}V_{k+1}$
are orthonormal as well as~\eqref{eq:kap}, we obtain
$$
\|V_{k+1} t_k\|_2
  \leq \|W^{-\sfrac{1}{2}}\|_2 \|W^{\sfrac{1}{2}}V_{k+1}t_k\|_2
  =    \|W\|_2^{-\sfrac{1}{2}} \|t_k\|_2
  \leq \sqrt{1+\delta} \|t_k\|_2.
$$
In GMRES, $\|t_k\|_2\rightarrow 0$ with increasing $k$,
which implies that $\|V_{k+1} t_k\|_2\rightarrow 0$ as well.
Therefore, we focus on bounding the residual gap $\|r_k - V_{k+1}t_k\|_2$
in order to satisfy~\eqref{eq:gap1} and~\eqref{eq:gap2}.

Suppose the matrix-vector products in the Arnoldi algorithm
are computed inexactly, i.e., line~4 in Algorithm~\ref{alg:arnoldi} 
is replaced by
\be\label{eq:epsilon}
w_j = (A+\cE_j)v_j, 
\ee
where $\|\cE_j\|_2\leq \epsilon_j$ for some given tolerance $\epsilon_j$.
Then in~\eqref{eq:3},
\be\label{eq:Ek}
E_k = \bmx \cE_1v_1, \cE_2v_2, \dots, \cE_kv_k \emx.
\ee
The following proposition bounds the residual gap at step $k$
in terms of the tolerances  $\epsilon_j$, for $j=1,\dots,k$. This is a direct corollary of  
results in \cite{SimoSzyl03} and~\cite{VDEhSlei04}.

\medskip
\begin{proposition}\label{thm:rgap}
Suppose that the inexact Arnoldi relation~\eqref{eq:3} holds,
where $E_k$ is given in~\eqref{eq:Ek}
with $\|\cE_j\|_2\leq\epsilon_j$ for $j=1,\dots,k$.
Then the resulting residual gap satisfies
\be\label{eq:rgap}
\|r_k-V_{k+1}t_k\|_2 \leq \|H_k^\dag\|_2 \sum_{j=1}^k \epsilon_j\|t_{j-1}\|_2.
\ee
\end{proposition}

\medskip

\subsection{A strategy for picking the $\epsilon_{j}$}

Proposition~\ref{thm:rgap} suggests the following strategy for picking
the tolerances $\epsilon_j$ that bound
the level of inexactness $\|{\cal E}_j\|_2$
in the matrix-vector products in~\eqref{eq:epsilon}.
Similarly to Theorem~\ref{thm:picketa}, let~$\phi_j$
be any positive numbers such that $\sum_{j=1}^k \phi_j = 1$.
If for all steps $j=1,\dots,k$, 
\be\label{eq:pickepsilon}
\epsilon_j \leq \frac{\phi_j\epsilon\sigma_{\min}(H_k)}{2}
                \frac{\|b\|_2}{\|t_{j-1}\|_2}, 
\ee
then from~\eqref{eq:rgap}
the residual gap in~\eqref{eq:gap1} satisfies
$$
\|r_k-V_{k+1}t_k\|_2 \leq \frac{\epsilon}{2} \|b\|_2.
$$
Interestingly, this result is independent of the accuracy of the inner products.
Similarly to~\eqref{eq:picketa}, the criterion for picking $\epsilon_j$ at step~$j$
involves $H_k$ that is only available at the final step~$k$.
A large number of numerical experiments~\cite{VDEhSlei04,BourFray06}
indicate that $\sigma_{\min}(H_k)$ can often be replaced by $1$.
Absorbing the factor $\sfrac{\phi_j}{2}$
into $\epsilon$ in~\eqref{eq:pickepsilon}
and replacing $\sigma_{\min}(H_k)$ by~$1$ or by~$\sigma_{\min}(A)$ 
leads, respectively, to the same aggressive and conservative thresholds
for $\epsilon_j$ as we obtained for $\eta_j$
in~\eqref{eq:picketa2} and in~\eqref{eq:picketa3}.
This suggests that matrix-vector products and inner products in GMRES
can be computed with the same level of inexactness.
We illustrate this with numerical examples in the next section.

\section{Numerical examples with emulated accuracy}\label{sec:num}

We illustrate our results with a few numerical examples.
We run GMRES with different matrices $A$ and right-hand sides $b$,
and compute the inner products and matrix-vector products inexactly
as in~\eqref{eq:eta} and~\eqref{eq:epsilon},
as described in Algorithm~\ref{alg:inexGMRES} below.

Note that the inner product $h_{j+1,j}$
in line 17 of Algorithm~\ref{alg:inexGMRES}
is also computed inexactly.
In Section~\ref{sec:IP}, to simplify the analysis,
we supposed that each $v_{j+1}$ was normalized exactly.
However, our numerical experiments indicate that $h_{j+1,j}$
can be computed with the same level of inexactness 
as the other inner products at step $j$.

We pick $\eta_{ij}$ randomly, 
uniformly distributed between $-\eta_j$ and $\eta_j$, 
and pick ${\cal E}_{j}$ to be a matrix of independent standard 
normal random variables, scaled to have norm~$\epsilon_{j}$.
Thus we have
$$
|\eta_{ij}|\leq \eta_{j}, \qquad
\|{\cal E}_{j}\|_2\leq\epsilon_{j},
$$
for chosen tolerances $\eta_j$ and $\epsilon_j$.
Throughout this first set of experiments, we use the same level of inexactness
for inner products and matrix-vector products,
i.e., $\eta_j=\epsilon_j$.

In the associated figures, the solid curve
is the relative residual $\|b-Ax_k\|_2/\|b\|_2$.
For reference, the dashed curve
is the relative residual if GMRES is run in double precision.
The crossed curve corresponds to 
the loss of orthogonality $\|F_k\|_2$ in~\eqref{eq:IF}.
The dotted curve is the chosen tolerance $\eta_j$.

\begin{algorithm}
\caption{A variable precision GMRES}
\label{alg:inexGMRES}
\begin{algorithmic}[1]
\REQUIRE{$A\in\R^{n\times n}$, $b\in\R^n$,
    $\epsilon>0$, $K_{\max}\in\mathbb{N}$, Conservative $\in\{0,1\}$}
\IF{Conservative} 
\STATE Compute or estimate $\sigma_{\min}(A)$
\ENDIF
\STATE $\beta = \sqrt{  b^Tb }$
\STATE $v_1 = b/\beta$
\FOR{$j=1,2,\dots,K_{\max}$}
\IF{Conservative} 
\STATE Compute $\eta_j$ and $\epsilon_j$ according to the bound (\ref{eq:picketa3})
\ELSE
\STATE Compute $\eta_j$ and $\epsilon_j$ according to the bound (\ref{eq:picketa2})
\ENDIF
\STATE Compute $w_{j} = (A+{\cal E}_j)v_j$ with $\Vert {\cal E}_j\Vert_2 \leq \epsilon_j $
\FOR{$i=1,\dots,j$}
\STATE Compute $h_{ij} =  v_i^T w_j +\eta_{i,j}$ with $\vert \eta_{i,j} \vert \leq \eta_j$
\STATE $w_{j} = w_{j} - h_{ij}v_i$
\ENDFOR
\STATE Compute $h_{j+1,j} = \sqrt{  w_{j}^T w_{j} + \eta_{j+1,j} }$
           with $\vert \eta_{j+1,j} \vert \leq \eta_j$
\IF{$h_{j+1,j}=0$}
\STATE Goto 27
\ENDIF
\STATE $v_{j+1} = w_{j} / h_{j+1,j}$
\STATE Compute $y_j$ and $\|t_j\|_2$, 
       the solution and residual of $\min_{y\in\R^j} \|\beta e_1 - H_j y\|_2$. 
\IF{$\|t_j\|_2<\epsilon$}
\STATE Goto 27
\ENDIF
\ENDFOR
\STATE Set $x_j=V_j y_j$
\end{algorithmic}
\end{algorithm}

\subsection{Relationship between $\eta_{ij}$ and loss of orthogonality}

Our first example illustrates the relationship between
the errors $\eta_{ij}$ in the inner products 
and the loss of orthogonality in the GMRES algorithm.

\begin{figure}[ht!] 
\begin{center}
\begin{tabular}{cc}
\includegraphics[width=2.3in]{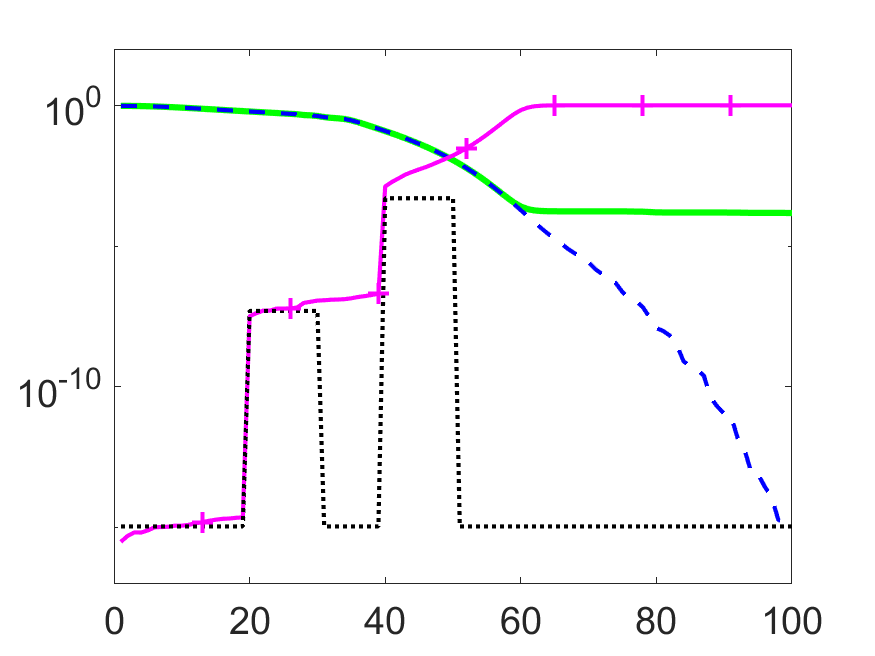} &
\includegraphics[width=2.3in]{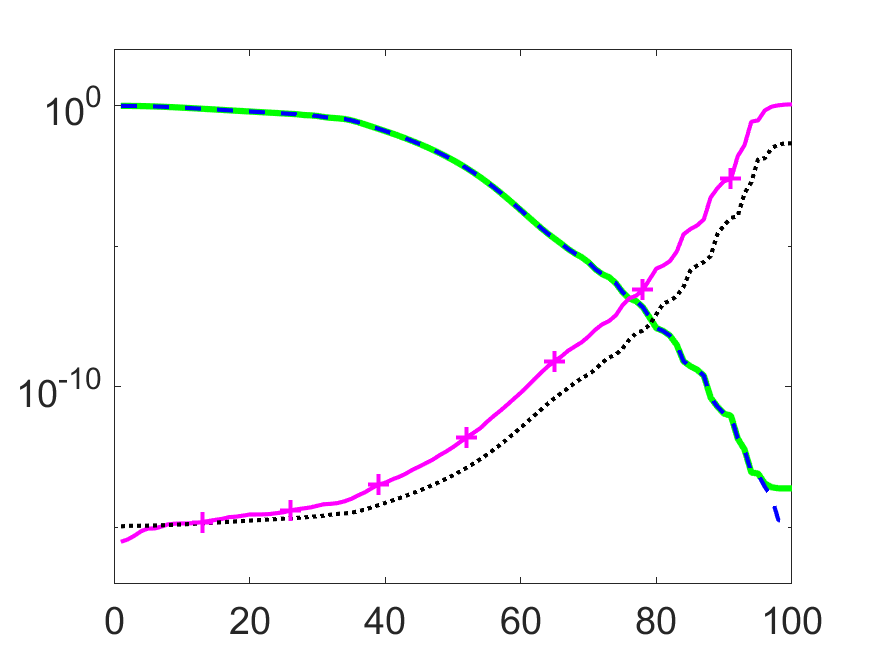} \\
{\footnotesize Example 1(a)} & 
{\footnotesize Example 1(b)} \\~~\\
\end{tabular}
\includegraphics[width=2.2in]{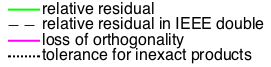} 
\caption{GMRES in variable precision: Grcar matrix.}
\label{fig:eg1}
\end{center}
\end{figure}

In this example, $A$ is the $100\times 100$ Grcar matrix of order $5$.
This is a highly non-normal Toeplitz matrix.
The right hand side is $b=A[\sin(1),\dots,\sin(100)]^T$.
Results are shown in Figure~\ref{fig:eg1}.

\noindent
In Example 1(a), 
$$
\eta_j = \epsilon_j =
\begin{cases}
 10^{-8}\|A\|_2, & \ \mbox{for} \ 20\leq j \leq 30, \\
 10^{-4}\|A\|_2, & \ \mbox{for} \ 40\leq j \leq 50, \\
 2^{-52}\|A\|_2, & \ \mbox{otherwise}.
\end{cases}
$$
The large increase in the inexactness of the inner products
at iterations $20$ and $40$ immediately 
leads to a large increase in~$\|F_k\|_2$.
This clearly illustrates the connection between
the inexactness of the inner products
and the loss of orthogonality in the Arnoldi vectors.
As proven in Theorem~\ref{thm:picketa},
until $\|F_k\|_2\approx 1$, the residual norm
is the same as it would have been
had all computations been performed in double precision.
Due to its large increases at iterations $20$ and $40$,
$\|F_k\|_2$ approaches $1$, and the residual norm starts to stagnate,
long before the relative residual norm reaches
the double precision machine precision.

In Example 1(b), the tolerances are chosen 
according to the aggressive criterion~\eqref{eq:picketa2} 
with $\epsilon=2^{-52}\|A\|_2$.
With this choice, $\|F_k\|_2$ does not reach $1$, 
and the residual norm does not stagnate until convergence.

\subsection{Conservative vs aggressive thresholds}

In our second example, $A$ is the matrix $\mathrm{494\_bus}$
from the SuiteSparse matrix collection~\cite{DaviHu11}.
This is a $494\times 494$ matrix with condition number $\kappa_2(A)\approx10^6$.
The right hand side is once again $b=A[\sin(1),\dots,\sin(100)]^T$.

\begin{figure}[ht!] 
\begin{center}
\begin{tabular}{cc}
\includegraphics[width=2.3in]{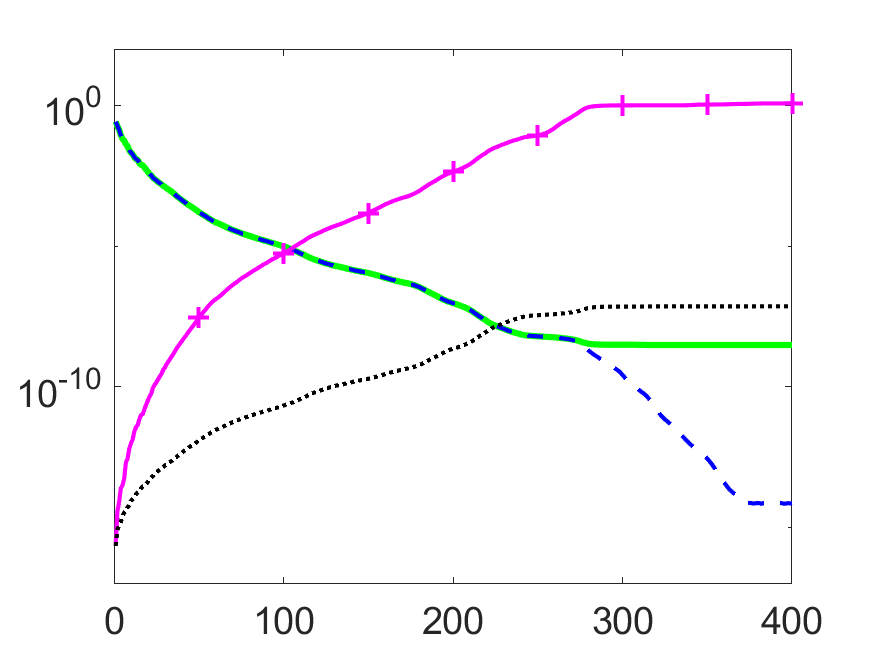} &
\includegraphics[width=2.3in]{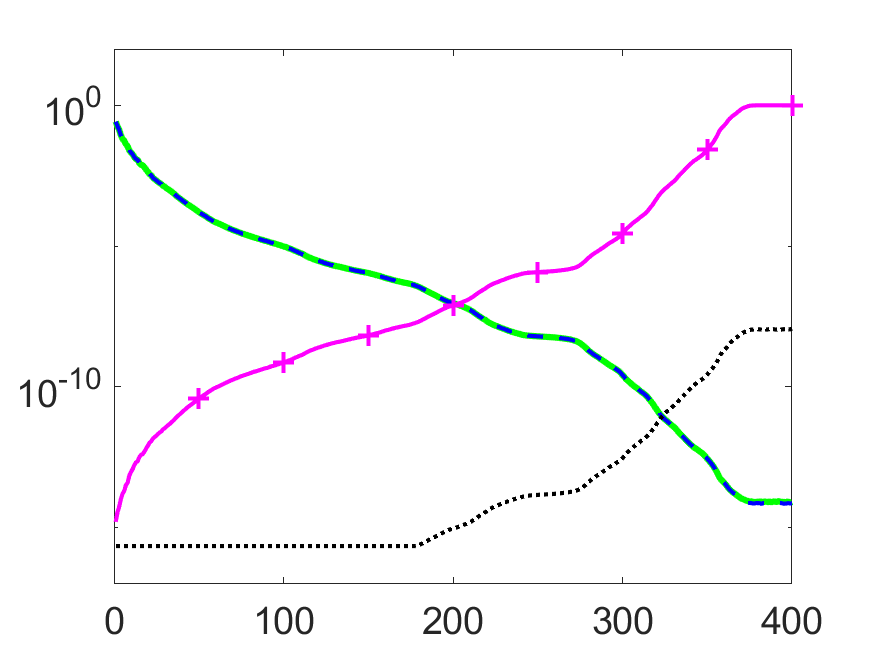} \\
{\footnotesize Example 2(a)} & 
{\footnotesize Example 2(b)} \\~~\\
\end{tabular}
\includegraphics[width=2.2in]{./legend.png} 
\caption{GMRES in variable precision: $\mathrm{494\_bus}$ matrix}
\label{fig:eg2}
\end{center}
\end{figure}

Results are shown in Figure~\ref{fig:eg2}.
In Example 2(a), tolerances are chosen according to
the aggressive threshhold~\eqref{eq:picketa2}
with $\epsilon=2^{-52}\|A\|_2$.
In this more ill-conditioned problem,
the residual norm starts to stagnate before
full convergence to double precision.
In Example 2(b), the tolerances are chosen according to
the conservative threshhold~\eqref{eq:picketa3}
with $\epsilon=2^{-52}\|A\|_2$, and there is no more such stagnation.
Because of these lower tolerances,
the inner products and matrix-vector products
have to be performed in double precision until about iteration 200.
This example illustrates the tradeoff between 
the level of inexactness and the maximum attainable accuracy.
The more ill-conditioned the matrix $A$ is,
the less opportunity there is for performing floating-point
operations inexactly in GMRES.

\section {Numerical experiments using variable floating-point arithmetic}\label{sec:mult}

In this section, we run variable floating-point arithmetic in order to assess the performances of
the approach in the context of half, simple and double precisions. These experiments are done with the Julia\footnote{https://julialang.org/}
language which allows to switch on demand between different floating-point types (Float16, Float32 and Float64).
We compute the inner products and matrix-vector products in  lower floating-point precision $u^{(low)}$
once either the aggressive threshold in~\eqref{eq:picketa2} or the conservative threshold in~\eqref{eq:picketa3} 
has increased above $u^{(low)}\|A\|_2$. 

\begin{figure}
\begin{center}
\begin{tabular}{rcc}
 &Conservative threshold~\eqref{eq:picketa3}&Aggressive threshold~\eqref{eq:picketa2}\\
 &&\\
\rotatebox{90}{ \hspace{0.5cm}Matrix-vector products}&\includegraphics[width=0.44\textwidth,keepaspectratio=true]{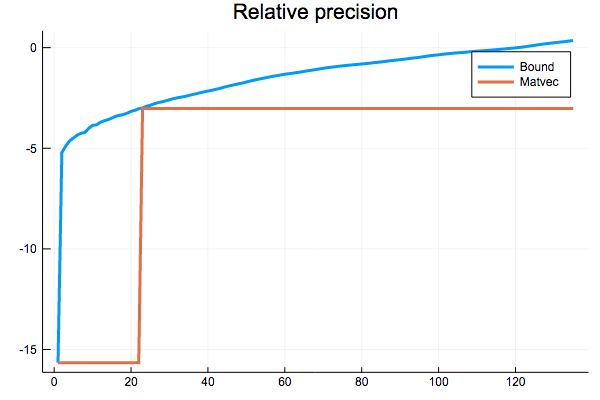}&
\includegraphics[width=0.44\textwidth,keepaspectratio=true]{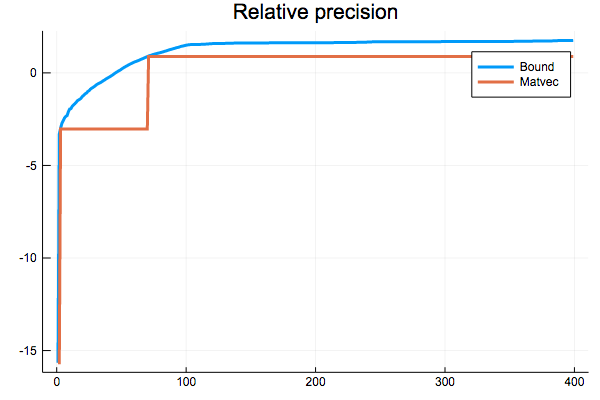}\\
 &&\\
\rotatebox{90}{  \hspace{1cm}Inner products}&\includegraphics[width=0.44\textwidth,keepaspectratio=true]{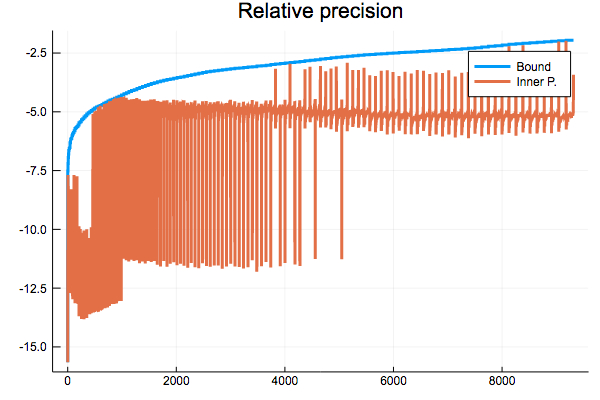}&\includegraphics[width=0.44\textwidth,keepaspectratio=true]{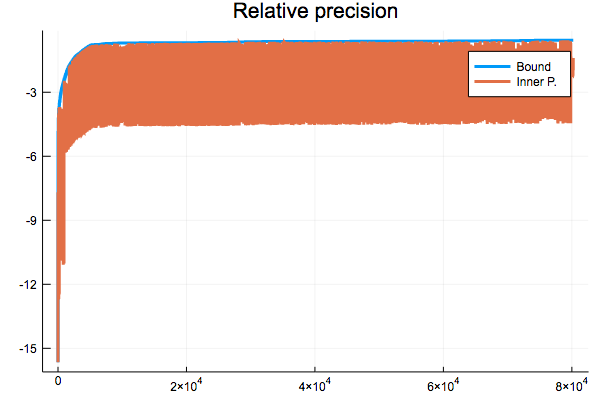}\\
 &&\\
\rotatebox{90}{  \hspace{0.75cm}Loss of orthogonality}&\includegraphics[width=0.44\textwidth,keepaspectratio=true]{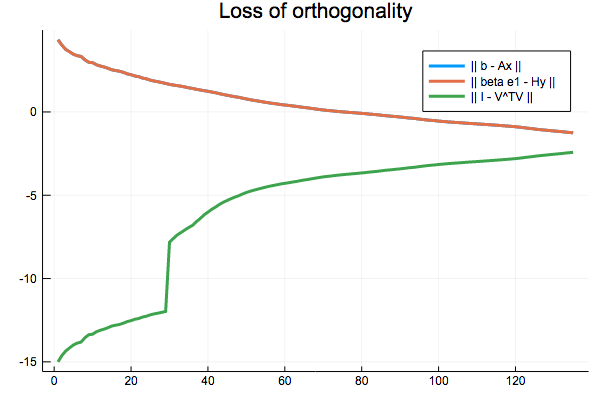}&\includegraphics[width=0.44\textwidth,keepaspectratio=true]{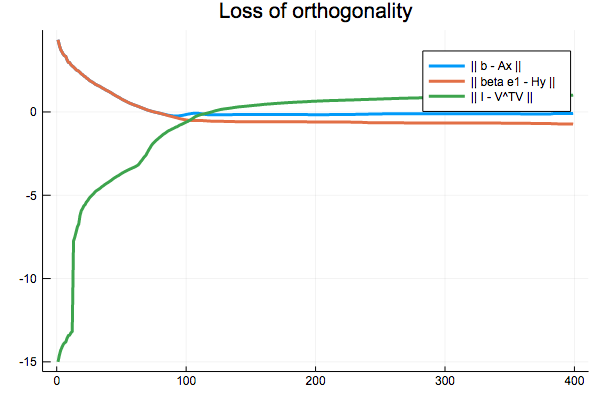}\\
\end{tabular}
\end{center}
\begin{center}
\begin{tabular}{cc}
\rotatebox{90}{  \hspace{0.75cm}Double precision}&\includegraphics[width=0.44\textwidth,keepaspectratio=true]{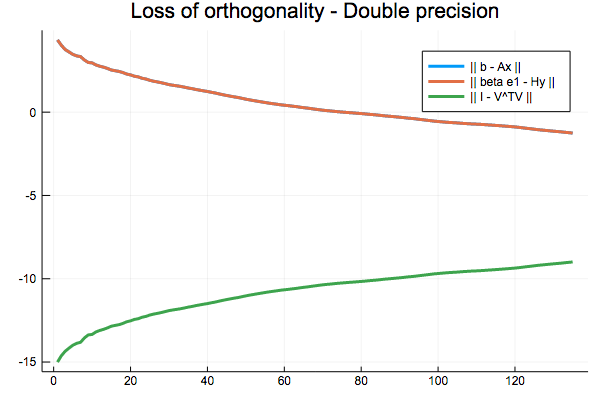}\\
\end{tabular}
\end{center}
 \caption{GMRES with variable precision floating-point arithmetic: Experiments with the $\mathrm{494\_bus}$ matrix and $\epsilon=10^{-6}\Vert A\Vert_2$. The bottom figure corresponds to GMRES in double precision.}
\label{fig:multi1} 
\end{figure}

\begin{figure}
\begin{center}
\begin{tabular}{rcc}
 &Conservative threshold~\eqref{eq:picketa3}&Aggressive threshold~\eqref{eq:picketa2}\\
 &&\\
\rotatebox{90}{ \hspace{0.5cm}Matrix-vector products}&\includegraphics[width=0.44\textwidth,keepaspectratio=true]{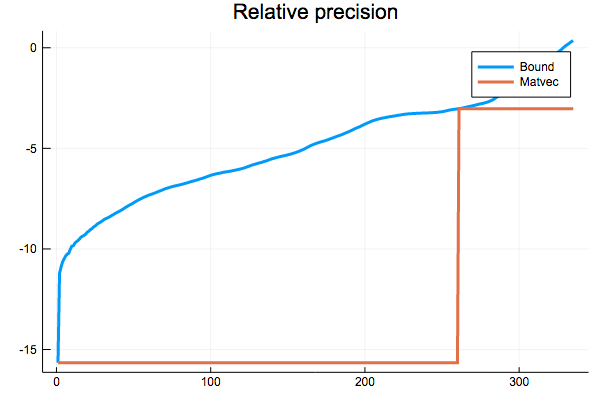}&
\includegraphics[width=0.44\textwidth,keepaspectratio=true]{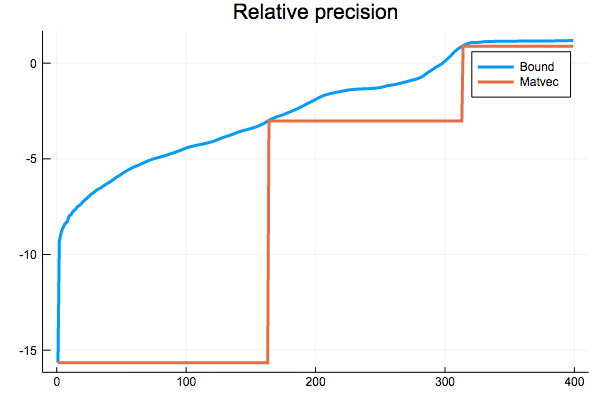}\\
 &&\\
\rotatebox{90}{ \hspace{1cm} Inner products}&\includegraphics[width=0.44\textwidth,keepaspectratio=true]{494bus-conservative-uphi-norm-eps6-innerP.png}&\includegraphics[width=0.44\textwidth,keepaspectratio=true]{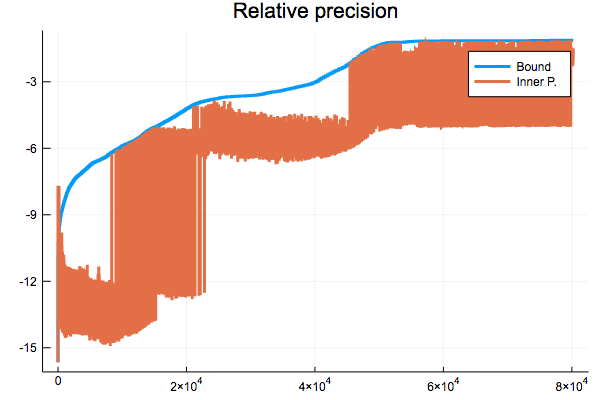}\\
 &&\\
\rotatebox{90}{ \hspace{0.75cm}Loss of orthogonality}&\includegraphics[width=0.44\textwidth,keepaspectratio=true]{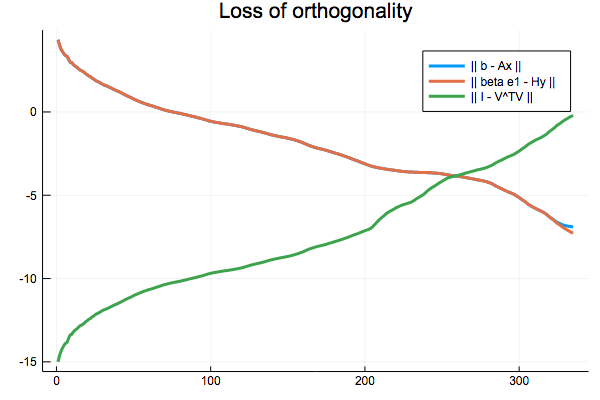}&\includegraphics[width=0.44\textwidth,keepaspectratio=true]{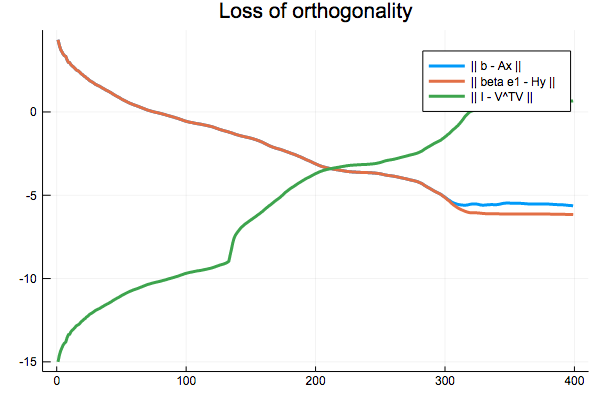}\\
\end{tabular}
\end{center}
\begin{center}
\begin{tabular}{cc}
\rotatebox{90}{  \hspace{0.75cm}Double precision}&\includegraphics[width=0.44\textwidth,keepaspectratio=true]{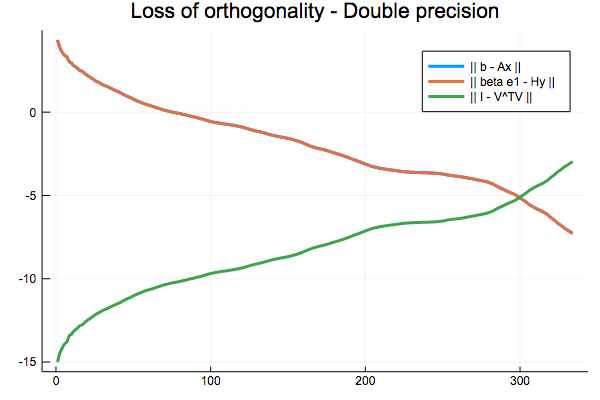}\\
\end{tabular}
\end{center}
 \caption{GMRES with variable precision floating-point arithmetic: Experiments with the $\mathrm{494\_bus}$ matrix and $\epsilon=10^{-12}\Vert A\Vert_2$. The bottom figure corresponds to GMRES in double precision.}
 \label{fig:multi2}
\end{figure}

Compared to the previous experiments, we must now take into account the magnitudes of both the vectors  and matrix perturbations  in order to select the precision of the computation. When no floating-point overflow arise, the perturbation of the matrix is simply computed from the difference between the norm of the matrix stored in Float64 and its conversion to Float16 and Float32. Regarding the inner products, we estimate their magnitude based on the sum of the exponents of the two vectors involved (plus 1 for the product of the mantissa). This explains the oscillating behavior of the accuracies observed in Figures~\ref{fig:multi1} and \ref{fig:multi2}: even if the precision has not changed, the estimated amplitude of the value of the inner products induces changes in the associated perturbations \eqref{eq:picketa2} and \eqref{eq:picketa3}. 

We focus only on the $494\_\mbox{bus}$ matrix using both the conservative and agressive thresholds.  The tolerances are chosen 
equal to $\epsilon=10^{-6}\|A\|$ and $\epsilon=10^{-12}\|A\|$ in order to illustrate the potential of the algorithm when moderate and high accuracies are required. In Figure~\ref{fig:multi1}, when a moderate decrease of the internally-recurred residual is required, we note that the conservative threshold results in a quick degradation of the precision for both the matrix-vector and inner products. All the matrix-vector products are computed in simple precision after 20 iterations, while the inner products start to be computed in simple precision after 30 iterations, precision that is mostly used after 90 iterations. We note a jump in the loss of orthogonality when the simple precision is triggered in the computation of the inner products. However, as expected from the theory, this does not degrade the decrease of the residual which is similar to the one observed with GMRES in double precision. When the aggressive threshold is used, we note that the simple precision is triggered after a few iterations for both the matrix-vector and inner products. The matrix-vector products are then computed in half precision from iteration 90 to convergence, while the precision of the inner products oscillate between  half and single depending on the amplitude of the vectors. The consequences are a complete loss of orthogonality after 100 iterations, which results in a slowing down in the decrease of the internally-recurred residual and a stagnation of the residual.

The effect of requiring a higher accuracy is mainly a delay in the exploitation of the multi-arithmetic. The results shown in Figure~\ref{fig:multi2} are similar to those obtained with a coarser tolerance, except for the delay in triggering the computation in simple, and half precisions.  We note again a similar decrease in the residuals compared to the GMRES algorithm in double precision. The loss of orthogonality is more severe when the agressive threshold is used, due to an earlier use of the simple precision in the inner products, as well as the use of the half precision in the latest iterations. The residual does not decrease anymore until the maximum number of iterations is reached.

This once again illustrates the tradeoff between 
the level of inexactness of the computations
and the maximum attainable accuracy.

\section{Conclusion}

We have shown how inner products
can be performed inexactly in MGS-GMRES without affecting
the convergence or final achievable accuracy of the algorithm.
We have also shown that a known framework
for inexact matrix-vector products is still valid
despite the loss of orthogonality in the Arnoldi vectors.
It would be interesting to investigate the impact
of scaling or preconditioning on these results.
Additionally, in future work, we plan to address
the question of how much computational savings can be achieved
by this approach on massively parallel computer architectures.

\newpage
\appendix

\section{Appendix}

\subsection{Proof of \eqref{eq:NUR}}

In line~7 of Algorithm~\ref{alg:arnoldi}, 
in the $\ell$th pass of the inner loop at step $j$, we have 
\be\label{eq:w1}
w_j^{(\ell)} = w_j^{(\ell-1)} - h_{\ell j} v_\ell
\ee
for $\ell=1,\dots,j$ and with $w_j^{(0)}=Av_j$.
Writing this equation for $\ell=i+1$ to $j$, we have
\begin{align*}
w_j^{(i+1)} &= w_j^{(i)}   - h_{i+1,j} v_{i+1}, \\
w_j^{(i+2)} &= w_j^{(i+1)} - h_{i+2,j} v_{i+2}, \\
&\ \ \vdots \\
w_j^{(j)}   &= w_j^{(j-1)} - h_{j,j} v_{j}.
\end{align*}
Summing the above and cancelling identical terms
that appear on the left and right hand sides gives
$$
w_j^{(j)} = w_j^{(i)} - \sum_{\ell=i+1}^{j} h_{\ell j} v_\ell.
$$
Because $w_j^{(j)} = v_{j+1} h_{j+1,j}$, this reduces to
\be\label{eq:w2}
w_j^{(i)} = \sum_{\ell=i+1}^{j+1} h_{\ell j} v_\ell.
\ee
Because the inner products $h_{ij}$ are computed inexactly as in~\eqref{eq:eta},
from~\eqref{eq:w1} we have
\begin{align*}
w_j^{(i)}
  &= w_j^{(i-1)} - h_{ij} v_i \\
  &= w_j^{(i-1)} - (v_i^T w_j^{(i-1)} + \eta_{ij}) v_i \\
  &= (I-v_iv_i^T)w_j^{(i-1)} - \eta_{ij}v_i.  
\end{align*}
Therefore,
$$
v_i^T w_j^{(i)} = - \eta_{ij}.
$$
Multiplying~\eqref{eq:w2} on the left by $-v_i^T$ gives
\be\label{eq:etaIJ}
\eta_{ij} = - \sum_{\ell=i+1}^{j+1} h_{\ell j} (v_i^Tv_\ell),
\ee
which is the entry in position $(i,j)$ of the matrix equation
$$
\bmx \eta_{11} & \dots  & \eta_{1k} \\
               & \ddots & \vdots    \\
               &        & \eta_{kk} \emx 
 = - \bmx v_1^Tv_2 & \dots  & v_1^Tv_{k+1} \\
                   & \ddots & \vdots       \\
                   &        & v_k^Tv_{k+1} \emx 
     \bmx h_{21} & \dots  & h_{2k}    \\
                 & \ddots & \vdots    \\
                 &        & h_{k+1,k} \emx, 
$$
i.e.,~\eqref{eq:NUR}.

\subsection{Proof of Lemma~\ref{lem:STLS}}

Equation~\eqref{eq:ejyk} follows from
$$
e_j^T H_k^\dag 
 \underbrace{\bmx H_{j-1} & 0 \\
                  0       & 0 \emx}_{\in\R^{(k+1)\times k}}
\bmx y_{j-1} \\ 0_{k-j+1} \emx						
  = e_j^T H_k^\dag H_k \bmx y_{j-1} \\ 0 \emx	
  = e_j^T \bmx y_{j-1} \\ 0 \emx=  0,
$$
and thus
\begin{align*}
|e_j^T y_k|
   = |e_j^T H_k^\dag \beta_1 e_1| 
  &= \left| e_j^T H_k^\dag
   \left(\beta_1 e_1 - \bmx H_{j-1} & 0 \\ 0 & 0 \emx \bmx y_{j-1} \\ 0 \emx \right) \right| \\
  &= \left| e_j^T H_k^\dag \bmx \beta_1 e_1 - H_{j-1} y_{j-1} \\ 0 \emx \right| 
   \leq \|H_k^\dag\|_2 \|t_{j-1}\|_2. 
\end{align*}
As for~\eqref{eq:STLS}, for any $\gamma>0$,
the smallest singular value of the matrix $\bmx \beta\gamma e_1, H_kD_k^{-1} \emx$
is the scaled total least squares (STLS) distance~\cite{PaigStra02b}
for the estimation problem $H_k D_k^{-1} z \approx \beta e_1$.
As shown in~\cite{PaigStra02a}, it can be bounded by
the least squares distance
$$
\min_z \|\beta e_1 - H_kD_k^{-1}z\|_2
 = \|\beta e_1 - H_kD_k^{-1}z_k\|_2
 = \|\beta e_1 - H_ky_k\|_2
 = \|t_k\|_2,
$$
where $z_k=D_ky_k$.
From~\cite[Theorem 4.1]{PaigStra02a}, we have
\be\label{eq:STLS2}
\frac{\|t_k\|_2}{ \big( \gamma^{-2} + \|D_ky_k\|_2^2/(1-\tau_k^2) \big)^{1/2} }
  \leq \sigma_{\min} \left( \bmx \beta\gamma e_1, H_kD_k^{-1} \emx \right)
  \leq \gamma\|t_k\|_2,
\ee
provided $\tau_k<1$, where
$$
\tau_k \equiv
   \frac{ \sigma_{\min}\left( \bmx \beta\gamma e_1, H_kD_k^{-1} \emx \right) }
        { \sigma_{\min} \left( H_kD_k^{-1} \right) }. 
$$
We now show that if $\gamma=(\epsilon\|b\|_2)^{-1}$
and $D_k$ satisfies~\eqref{eq:condD}, then $\tau_k\leq\sfrac{1}{\sqrt{2}}$.
From the upper bound in~\eqref{eq:STLS2} we immediately have
$$
\sigma_{\min}\left( \bmx \beta\gamma e_1, H_kD_k^{-1} \emx \right)
   \leq \gamma\|t_k\|_2 = \frac{\|t_k\|_2}{\epsilon\|b\|_2}.
$$
Also,
$$
\sigma_{\min}\left( H_kD_k^{-1} \right)
   = \min_{z\neq 0} \frac{\|H_kD_k^{-1} z\|_2}{\|z\|_2}
   = \min_{z\neq 0} \frac{\|H_kz\|_2}{\|D_kz\|_2}
   \geq \min_{z\neq 0} \frac{\|H_kz\|_2}{\|D_k\|_2\|z\|_2}
   = \frac{\sigma_{\min}(H_k)}{\|D_k\|_2}.
$$
Therefore, if~\eqref{eq:condD} holds,
$$
\tau_k \leq \frac{\|t_k\|_2}{\epsilon\|b\|_2}\frac{\|D_k\|_2}{\sigma_{\min}(H_k)} 
       \leq \frac{1}{\sqrt{2}}.
$$
Substituting $\gamma=(\epsilon\|b\|_2)^{-1}$ and $\tau_k\leq\sfrac{1}{\sqrt{2}}$
into~\eqref{eq:STLS2} gives~\eqref{eq:STLS}.

\subsection{Proof of Lemma~\ref{lem:qre}}

Note from~\eqref{eq:q1} that the singular values
of $Q$ satisfy
$$
\big( \sigma_i(Q) \big)^2
  = \sigma_i(Q^TQ) = \sigma_i( I_k - F ),
  \qquad i=1,\dots,k.
$$
Therefore,
\be\label{eq:QF}
\sqrt{1 - \|F\|_2} \leq \sigma_i(Q) 
  \leq \sqrt{1 + \|F\|_2}, 
  \qquad i=1,\dots,k.
\ee
Equation~\eqref{eq:q2} is equivalent
to the linear matrix equation
$$
Q^TMQ = I_k - Q^TQ.
$$
It is straightforward to verify that
one matrix~$M$ satisfying this equation is 
\begin{align*}
M &= (Q^\dag)^T (I_k - Q^TQ) Q^\dag \\
  &= Q(Q^TQ)^{-1} (I_k - Q^TQ) (Q^TQ)^{-1} Q^T.
\end{align*}
Notice that the above matrix $M$ is symmetric.
It can also be verified
using the singular value decomposition of $Q$
that the eigenvalues and singular values of $I_n+M$ are
$$
\lambda_i(I_n+M)
   = \sigma_i(I_n+M) 
   = \begin{cases}
     \big( \sigma_i(Q) \big)^{-2},  & \quad i=1,\dots,k, \\
     1,                             & \quad i=k+1,\dots,n, 
     \end{cases}
$$
which implies that the matrix $I_n+M$ is positive definite.
From the above and~\eqref{eq:QF}, provided $\|F\|_2\leq\delta<1$,
\be\label{eq:kap}
\frac{1}{1+\delta}
  \leq \frac{1}{\big( \sigma_{\max}(Q) \big)^2}
  \leq \sigma_{i}(I_n+M)
  \leq \frac{1}{\big( \sigma_{\min}(Q) \big)^2}
  \leq \frac{1}{1-\delta},
\ee
from which~\eqref{eq:kapdelta} follows.

\section*{Acknowledgments}

The authors would like to thank two anonymous referees whose
comments lead to significant improvements in the presentation.

\bibliographystyle{siamplain}
\bibliography{refs}

\begin{thebibliography}{10}

\bibitem{Bjor67}
{\sc A.~Bj{\"o}rck}, {\em Solving linear least squares problems by
  {G}ram-{S}chmidt orthogonalization}, BIT Numerical Mathematics, 7 (1967),
  pp.~1--21.

\bibitem{BjorPaig92}
{\sc A.~Bj{\"o}rck and C.~Paige}, {\em Loss and recapture of orthogonality in
  the {M}odified {G}ram-{S}chmidt algorithm}, SIAM Journal on Matrix Analysis
  and Applications, 13 (1992), pp.~176--190.

\bibitem{BourFray06}
{\sc A.~Bouras and V.~Frayss{\'e}}, {\em Inexact matrix-vector products in
  {K}rylov methods for solving linear systems: A relaxation strategy}, SIAM
  Journal on Matrix Analysis and Applications, 26 (2006), pp.~660--678.

\bibitem{CarsHigh18}
{\sc E.~Carson and N.J.Higham}, {\em Accelerating the solution of linear
  systems by iterative refinement in three precisions}, SIAM Journal on
  Scientific Computing, 40 (2018), pp.~817--847.

\bibitem{DaviHu11}
{\sc T.~A. Davis and Y.~Hu}, {\em The {U}niversity of {F}lorida sparse matrix
  collection}, ACM Transactions on Mathematical Software, 38 (2011), pp.~1--25.

\bibitem{Drko95}
{\sc J.~Drko{\v s}ov{\'a}, A.~Greenbaum, M.~Rozlo{\v z}n{\' i}k, and
  Z.~Strako{\v s}}, {\em Numerical stability of {GMRES}}, BIT Numerical
  Mathematics, 35 (1995), pp.~309--330.

\bibitem{VDEhSlei04}
{\sc J.~V.~D. Eshof and G.~Sleijpen}, {\em Inexact {K}rylov subspace methods
  for linear systems}, SIAM Journal on Matrix Analysis and Applications, 26
  (2004), pp.~125--153.

\bibitem{Fre92b}
{\sc R.~Freund}, {\em Quasi-kernel polynomials and convergence results for
  quasi-minimal residual iterations}, Numerical Methods in Approximation
  Theory, 9 (1992), pp.~77--95.

\bibitem{Fre92a}
{\sc R.~Freund}, {\em Quasi-kernel polynomials and their use in non-{Hermitian}
  matrix iterations}, Journal of Computational and Applied Mathematics, 43
  (1992), pp.~135--158.

\bibitem{FreNac91}
{\sc R.~Freund and N.~Nachtigal}, {\em {QMR}: {A} quasi-minimal residual method
  for non-{Hermitian} linear systems}, Numerische Mathematik, 60 (1991),
  pp.~315--339.

\bibitem{GiraGratLang07}
{\sc L.~Giraud, S.~Gratton, and J.~Langou}, {\em Convergence in backward error
  of relaxed {GMRES}}, SIAM Journal on Scientific Computing, 29 (2007),
  pp.~710--728.

\bibitem{GratSimoToin18}
{\sc S.~Gratton, E.~Simon, and P.~Toint}, {\em Minimizing convex quadratic with
  variable precision {K}rylov methods}, ar{X}iv, abs/1807.07476 (2018).

\bibitem{Gree97}
{\sc A.~Greenbaum, M.~Rozlo{\v z}n{\' i}k, and Z.~Strako{\v s}}, {\em Numerical
  behaviour of the {M}odified {G}ram-{S}chmidt {GMRES} implementation}, BIT
  Numerical Mathematics, 37 (1997), pp.~706--719.

\bibitem{haid2018b}
{\sc A.~Haidar, A.~Abdelfattah, M.~Zounon, P.~Wu, S.~Pradesh, S.~Tomov, and
  J.~Dongarra}, {\em The design of fast and energy-efficient linear solvers: on
  the potential of half-precision arithmetic and iterative refinement
  techniques}, in Computational Science—ICCS 2018, Yong Shi, Haohuan Fu,
  Yingjie Tian, Valeria V. Krzhizhanovskaya, Michael Harold Lees, Jack
  Dongarra, and Peter M. A. Sloot editors, Springer, 2018, pp.~586--600.

\bibitem{haid2018}
{\sc A.~Haidar, S.~Tomov, J.~Dongarra, and N.~J. Higham}, {\em Harnessing {GPU}
  tensor cores for fast {FP16} arithmetic to speed up mixed-precision iterative
  refinement solvers}, in Proceedings of the International Conference for High
  Performance Computing, Networking, Storage, and Analysis, IEEE Press, 2018,
  pp.~47:1--47:11.

\bibitem{haid2017}
{\sc A.~Haidar, P.~Wu, S.~Tomov, and J.~Dongarra}, {\em Investigating half
  precision arithmetic to accelerate dense linear system solvers}, in
  Proceedings of the 8th Workshop on Latest Advances in Scalable Algorithms for
  Large-Scale Systems, ScalA 17, 2017, pp.~10:1--10:8.

\bibitem{High17}
{\sc N.~J. Higham}, {\em A multiprecision world}, SIAM News, 50 (2017).

\bibitem{HighPran19b}
{\sc N.~J. Higham and S.~Pradesh}, {\em Simulating low precision floating-point
  arithmetic}, SIAM Journal on Scientific Computing, 41 (2019), pp.~585--602.

\bibitem{HighPran19}
{\sc N.~J. Higham and S.~Pradesh}, {\em Squeezing a matrix into half precision,
  with an application to solving linear systems}, SIAM Journal on Scientific
  Computing, 41 (2019), pp.~2536--2551.

\bibitem{LeonBjorGand13}
{\sc S.~J. Leon, A.~Bj{\"o}rck, and W.~Gander}, {\em Gram-{S}chmidt
  orthogonalization: 100 years and more}, Numerical Linear Algebra with
  Applications, 20 (2013), pp.~492--532.

\bibitem{Paig09}
{\sc C.~Paige}, {\em A useful form of unitary matrix obtained from any sequence
  of unit $2$-norm $n$-vectors}, SIAM Journal on Matrix Analysis and
  Applications, 31 (2009), pp.~565--583.

\bibitem{PaigRozlStra06}
{\sc C.~Paige, M.~Rozlo{\v z}n{\' i}k, and Z.~Strako{\v s}}, {\em Modified
  {G}ram-{S}chmidt ({MGS}), least squares, and backward stability of
  {MGS-GMRES}}, SIAM Journal on Matrix Analysis and Applications, 28 (2006),
  pp.~264--284.

\bibitem{PaigStra02a}
{\sc C.~Paige and Z.~Strako{\v s}}, {\em Bounds for the least squares distance
  using scaled total least squares}, Numerische Mathematik, 91 (2002),
  pp.~93--115.

\bibitem{PaigStra02b}
{\sc C.~Paige and Z.~Strako{\v s}}, {\em Scaled total least squares
  fundamentals}, Numerische Mathematik, 91 (2002), pp.~117--146.

\bibitem{PaigWull14}
{\sc C.~Paige and W.~W\"ulling}, {\em Properties of a unitary matrix obtained
  from a sequence of normalized vectors}, SIAM Journal on Matrix Analysis and
  Applications, 35 (2014), pp.~526--545.

\bibitem{PestWath13}
{\sc J.~Pestana and A.~J. Wathen}, {\em On the choice of preconditioner for
  minimum residual methods for non-{H}ermitian matrices}, Journal of
  Computational and Applied Mathematics, 249 (2013), pp.~57--68.

\bibitem{SaadSchu86}
{\sc Y.~Saad and M.~H. Schultz}, {\em {GMRES}: A generalized minimal residual
  algorithm for solving nonsymmetric linear systems}, SIAM Journal on
  Scientific and Statistical Computing, 7 (1986), pp.~856--869.

\bibitem{SimoSzyl03}
{\sc V.~Simoncini and D.~Szyld}, {\em Theory of inexact {K}rylov subspace
  methods and applications to scientific computing}, SIAM Journal on Scientific
  Computing, 25 (2003), pp.~454--477.

\bibitem{SimoSzyl07}
{\sc V.~Simoncini and D.~Szyld}, {\em Recent computational developments in
  {K}rylov subspace methods for linear systems}, Numerical Linear Algebra with
  Applications, 14 (2007), pp.~1--59.

\end{thebibliography}
\end{document}